\documentclass[11pt]{article}
\usepackage{geometry}                
\usepackage{graphicx}
\usepackage{amssymb}
\usepackage{epstopdf}
\usepackage{amsmath,amsthm,amscd,amssymb}
\usepackage{latexsym}
\numberwithin{equation}{section}

\theoremstyle{plain}
\newtheorem{theorem}{Theorem}[section]

\theoremstyle{definition}

\theoremstyle{remark}
\newtheorem{remark}[theorem]{Remark}

\newtheorem{case[theorem]}{Case}

\title{Sums and products in finite fields: an integral geometric viewpoint}
\author{Derrick Hart and Alex Iosevich}

\begin{document}

\maketitle

\begin{abstract} We prove that if $A \subset {\Bbb F}_q$ is such that
$$|A|>q^{\frac{1}{2}+\frac{1}{2d}},$$ then
$${\Bbb F}_q^{*} \subset dA^2=A^2+\dots+A^2 \ \ d \  \text{times},$$ where
$$A^2=\{a \cdot a': a,a' \in A\},$$ and where ${\Bbb F}_q^{*}$ denotes the
multiplicative group of the finite field ${\Bbb F}_q$. In particular, we
cover ${\Bbb F}_q^{*}$ by $A^2+A^2$ if $|A|>q^{\frac{3}{4}}$.
Furthermore, we prove that if
$$|A| \ge C_{size}^{\frac{1}{d}}q^{\frac{1}{2}+\frac{1}{2(2d-1)}},$$ then
$$|dA^2| \ge q \cdot \frac{C^2_{size}}{C^2_{size}+1}.$$

Thus $dA^2$ contains a positive proportion of the elements of ${\Bbb F}_q$
under a considerably weaker size assumption.We use the geometry of
${\Bbb F}_q^d$, averages over hyper-planes and orthogonality
properties of character sums. In particular, we see that using
operators that are smoothing on $L^2$ in the Euclidean setting leads
to non-trivial arithmetic consequences in the context of finite
fields.
\end{abstract}

\tableofcontents

\section{Introducion}

A classical problem in additive number theory is to determine, given a
finite subset $A$ of a ring, whether both $2A=\{a+a':a,a' \in A\}$ and
$A^2=\{a \cdot a': a,a' \in A\}$ can be small in a suitable sense. A related
question, posed in a finite field ${\Bbb F}_q$ with $q$ elements, is how
large $A \subset {\Bbb F}_q$ need to be to assure that
$dA^2=A^2+A^2+\dots+A^2={\Bbb F}_q$. It is known (see e.g. \cite{GK06}) that
if $d=3$ and $q$ is prime, this conclusion is assured if $|A| \ge
Cq^{\frac{3}{4}}$, with a sufficiently large constant $C>0$. It is
reasonable to conjecture that if $|A| \ge
C_{\epsilon}q^{\frac{1}{2}+\epsilon}$, then $2A^2={\Bbb F}_q$. This result
cannot hold, especially in the setting of general finite fields if
$|A|=\sqrt{q}$ because $A$ may in fact be a subfield. See also \cite{BGK06},
\cite{C04}, \cite{G06}, \cite{G07}, \cite{HIS07}, \cite{KS07}, \cite{TV06},
\cite{V07} and the references contained therein on recent progress related
to this problem and its analogs.

For example, it is proved in \cite{G06} that
$$ 8X \cdot Y={\Bbb Z}_p, $$ for $p$ prime, provided that $|X||Y|>p$ and
either $Y=-Y$ or $Y \cap (-Y)=\emptyset$. In \cite{GK06} the author
prove that if $A$ is subgroup of ${\Bbb Z}_p^{*}$, and
$|A|>p^{\delta}$, $\delta>0$, then
$$ NA={\Bbb Z}_p$$ with
$$ N \leq C4^{\frac{1}{\delta}}.$$

The purpose of this paper is to use the geometry of ${\Bbb F}_q^d$,
where $q$ is not necessarily a prime number, to deduce a good lower bound on the size of $A$ that guarantees that $dA^2={\Bbb
F}_q$, with the possible exception of $0$.  Furthermore, it is shown that the lower bound on $A$ may be relaxed if one settles for a positive proportion of $\Bbb F_q$. Our main result is the following.

\begin{theorem} \label{kickass} Let $A \subset {\Bbb F}_q$, where ${\Bbb
F}_q$ is an arbitrary finite field with $q$ elements, such that
$|A|>q^{\frac{1}{2}+\frac{1}{2d}}$. Then
\begin{equation} \label{sex} {\Bbb F}_q^{*} \subset dA^2.\end{equation}

Suppose that
$$|A| \ge C^{\frac{1}{d}}_{size} q^{\frac{1}{2}+\frac{1}{2(2d-1)}}.$$

Then
\begin{equation} \label{sex2} |dA^2| \ge q \cdot
\frac{C^{2-\frac{1}{d}}_{size}}{C^{2-\frac{1}{d}}_{size}+1}. \end{equation}
\end{theorem}

In particular, if $d=2$,
$${\Bbb F}_q^{*} \subset A^2+A^2$$ if
$$|A|>q^{\frac{3}{4}},$$ and
$$ |A^2+A^2| \ge q \cdot \frac{C^2_{size}}{C^2_{size}+1}$$ if
$$ |A| \ge C_{size}^{\frac{1}{2}} q^{\frac{2}{3}}.$$

Also, Theorem \ref{kickass} gives an explicit bound for the conjecture
mentioned in \cite{GK06}, namely that if $|A| \ge
C_{\epsilon}q^{\frac{1}{2}+\epsilon}$, there exists $d=d(\epsilon)$ such
that $dA^2$ covers ${\Bbb F}_q$. In view of this, we restate Theorem
\ref{kickass} as follows.

\begin{theorem} \label{kickass2} Let $A \subset {\Bbb F}_q$, where ${\Bbb
F}_q$ is an arbitrary finite field with $q$ elements, such that
$$|A| \ge C_{\epsilon} q^{\frac{1}{2}+\epsilon},$$ for some $\epsilon>0$.
Then (\ref{sex}) holds for
$d=d(\epsilon)$ equal to the smallest integer greater than or equal to
$\frac{1}{2 \epsilon}$. Moreover, (\ref{sex2}) holds if $d$ is equal to the
smallest integer greater than or equal to $\frac{1}{2}+\frac{1}{4\epsilon}$.
\end{theorem}

Throughout the paper, $X \lesssim Y$ means that there exists a universal
constant $C$, indepedent of $q$, such that $X \leq CY$, and $X \approx Y$
means that $X \lesssim Y$ and $Y \lesssim X$. In the instances when the size
of the constant matters, this fact shall be mentioned explicitly.

\begin{remark} The reader can easily check that in Theorem \ref{kickass} and
Theorem \ref{kickass2}, $dA^2$ may be easily replaced by
$$A _1 \cdot B_1+\dots+A_d \cdot B_d, $$ provided that
$$\Pi_{j=1}^d |A_j||B_j| \ge Cq^{d+1}$$ with a sufficiently large constant
$C>0$.
\end{remark}

The proof of Theorem \ref{kickass} is based on the following geometric
observation that is interesting in its own right.
\begin{theorem} \label{main} Let $E \subset {\Bbb F}_q^d$ such that
$|E|>q^{\frac{d+1}{2}}$. Then
$$ {\Bbb F}_q^{*} \subset \{x \cdot y: x,y \in E\}.$$
\end{theorem}

To prove Theorem \ref{kickass} we shall need the following conditional
version of Theorem \ref{main}.
\begin{theorem} \label{main+} Let $E \subset {\Bbb F}_q^d$ such that
$$ |E \cap l_y| \leq C_{geom}q^{\frac{\alpha}{d}}$$ for some $0 \leq \alpha
\leq d$, for every
$y \in {\Bbb F}_q^d$, $y \not=(0, \dots, 0)$, where
$$ l_y=\{ty: t \in {\Bbb F}_q\}.$$

Suppose that
$$ |E| \ge C_{size}q^{\frac{d}{2}+\frac{\alpha}{2d}}.$$

Then
$$ |\{x \cdot y: x,y \in E\}| \ge q \cdot \frac{
C_{size}^2}{C_{size}^2+C_{geom}}.$$
\end{theorem}

\vskip.125in

\begin{remark} Theorem \ref{main+} has non-trivial applications to many
other problems in additive number theory and geometric combinatorics, such
as the Erd\H os distance problem, distribution of simplexes and others. We
study these problems systematically in \cite{HIR07}.  \end{remark}

\subsection{Integral geometric viewpoint} At the core of the proof of
Theorem \ref{main} and Theorem \ref{main+} is the
$L^2({\Bbb F}_q^d)$ estimate for the "rotating planes" operator
$$ {\cal R}_tf(x)=\sum_{x \cdot y=t} f(y).$$

In the Euclidean space, this operator is a classical example of a
phenomenon, thoroughly explored by Hormander, Phong, Stein and others (see
e.g. \cite{St93}) and the references contained therein) where an operator
that averages a function over a family of manifolds satisfies better than
trivial bounds on $L^2({\Bbb F}_q^d)$ provided that the family of manifolds
satisfies an appropriate curvature condition. It turns out that in the
finite field setting, the aforementioned operator, suitably interpreted,
satisfies analogous bounds which lead to interesting arithmetic
consequences.

In contrast, the authors of \cite{HIS07} took advantage of the $L^2({\Bbb
F}_q^d)$ mapping properties of the operator
$$ H_jf(x)=\sum_{y_1y_2=j} f(x-y),$$ and in \cite{IR07} the underlying
operator is
$$ A_tf(x)=\sum_{y_1^2+\dots+y_d^2=t} f(x-y),$$ though in neither paper was
this perspective made explicit. These examples suggest that systematic
theory of Fourier Integral Operator in the setting of vector spaces over
finite fields needs to be worked out and the authors shall take up this task
in a subsequent paper.

\vskip.125in

\subsection{Fourier analysis used in this paper} Let $f: {\Bbb F}_q^d \to
{\Bbb C}$. Let $\chi$ be a non-trivial additive character on ${\Bbb F}_q$.
Define the Fourier transform of $f$ by the formula
$$ \widehat{f}(m)=q^{-d} \sum_{x \in {\Bbb F}_q^d} \chi(-x \cdot m) f(x)$$
for $m \in {\Bbb F}_q^d$.

The formulas we shall need are the following:
$$ \sum_{t \in {\Bbb F}_q} \chi(-at)=0 \ \ \text{(orthogonality)},$$ if $t
\not=0$, and $q$ otherwise,
$$ f(x)=\sum_m \chi(x \cdot m) \widehat{f}(m) \ \ \text{(inversion)},$$
$$ \sum_m \widehat{f}(m) \overline{\widehat{g}(m)}=q^{-d} \sum_x f(x)g(x)
\ \ \text{(Plancherel/Parseval)}.$$

In the case when $q$ is a prime, one may take $\chi(t)=e^{\frac{2 \pi
i}{q}t}$, and in the general case the formula is only slightly more
complicated.

\vskip.125in

\subsection{Acknowledgements:} The authors wish to thank Moubariz Garaev,
Nets Katz, Sergei Konyagin and Ignacio Uriarte-Tuero for a thorough
proofreading of the earlier drafts of this paper and for many interesting
and helpful remarks.

\vskip.125in

\section{Proof of the basic geometric estimate (Theorem \ref{main})}

\vskip.125in

Let
$$ \nu(t)=|\{(x,y) \in E \times E: x \cdot y=t\}|.$$

We have
$$ \nu(t)=\sum_{x,y \in E} q^{-1} \sum_{s \in {\Bbb F}_q} \chi(s(x \cdot
y-t)),$$ where $\chi$ is a non-trivial additive character on ${\Bbb F}_q$.
It follows that
$$ \nu(t)={|E|}^2q^{-1}+R, $$ where
$$ R=\sum_{x,y \in E} q^{-1} \sum_{s \not=0} \chi(s(x \cdot y-t)).$$

Viewing $R$ as a sum in $x$, applying the Cauchy-Schwartz inequality and
dominating the sum over $x \in E$ by the sum over $x \in {\Bbb F}_q^d$, we
see that
$$R^2 \leq |E| \sum_{x \in {\Bbb F}_q^d} q^{-2} \sum_{s,s' \not=0}
\sum_{y,y' \in E}
\chi(sx \cdot y-s'x \cdot y') \chi(t(s'-s)).$$

Orthogonality in the $x$ variable yields
$$=|E| q^{d-2} \sum_{\substack{sy=s'y' \\ s,s' \not=0}}
\chi(t(s'-s))E(y)E(y').$$

If $s \not=s'$ we may set $a=s/s', b=s'$ and obtain
$$ |E| q^{d-2} \sum_{\substack{y \not=y' \\ ay=y' \\ a \not=1,b}}
\chi(tb(1-a))E(y)E(y')$$
$$=-|E| q^{d-2} \sum_{y \not=y', a \not=1} E(y)E(ay), $$ and the absolute
value of this quantity is
$$ \leq |E|q^{d-2} \sum_{y \in E} |E \cap l_y|$$
$$ \leq {|E|}^2 q^{d-1},$$ since
$$|E \cap l_y| \leq q$$ by the virtue of the fact that each line contains
exactly $q$ points.

If $s=s'$ we get
$$ |E|q^{d-2} \sum_{s,y} E(y)={|E|}^2q^{d-1}.$$

It follows that
$$ \nu(t)={|E|}^2q^{-1}+R(t),$$ where
$$ R^2(t) \leq -Q(t)+{|E|}^2q^{d-1},$$ with
$$ Q(t) \ge 0.$$

It follows that
$$ R^2(t) \leq {|E|}^2q^{d-1},$$ so
\begin{equation} \label{remainder} |R(t)| \leq
|E|q^{\frac{d-1}{2}}.\end{equation}

We conclude that
$$ \nu(t)={|E|}^2q^{-1}+R(t)$$ with $|R(t)|$ bounded as in
(\ref{remainder}).

This quantity is strictly positive if $|E|>q^{\frac{d+1}{2}}$ with a
sufficiently large constant $C>0$. This completes the proof of Theorem
\ref{main}. Theorem \ref{kickass} follows from Theorem \ref{main} by simply
setting $E=A \times A \times \dots \times A$.

\vskip.125in

\section{Proof of the enhanced geometric estimate (Theorem \ref{main+})}

\vskip.125in

Assume throughout the argument, without loss of generality, that $E$ does not
contain the origin. Applying Cauchy-Schwartz as above we see that
$$ \nu^2(t) \leq |E| \sum_{x \in E} \sum_{y,y' \in E} q^{-2} \sum_{s,s'}
\chi(x \cdot (sy-s'y')) \chi(t(s'-s)).$$

It follows that
$$ \sum_t \nu^2(t) \leq |E|q^{d-1} \sum_{s} \sum_m \widehat{E}(sm)
\sum_{y-y'=m} E(y)E(y')$$
$$=|E|q^{d-1} \sum_{s} \sum_m \widehat{E}(ms) E*E(m)$$
\begin{equation} \label{stage1}=|E|q^{d-1} \sum_m  \left( \sum_{s}
\widehat{E}(sm)\right)
E*E(m).\end{equation}

Now,
$$ \sum_{s} \widehat{E}(ms)=\sum_s q^{-d} \sum_x E(x) \chi(-x \cdot ms)$$
$$=q^{-(d-1)} \sum_{x \cdot m=0} E(x).$$

Inserting this it into (\ref{stage1}) we get
\begin{equation} \label{stage2} |E| \sum_m \left( \sum_{x \cdot m=0} E(x)
\right) \cdot E*E(m). \end{equation}

Let
$$ F(m)=\sum_{x \cdot m=0} E(x), \ \ G(m)=E*E(m).$$

By a direct calculation,
$$ \widehat{G}(k)=q^d {|\widehat{E}(k)|}^2.$$

On the other hand,
$$ \widehat{F}(k)=q^{-d} \sum_m \chi(-m \cdot k) \sum_{x \cdot m=0} E(x)$$
$$=q^{-d}q^{-1} \sum_{m,x} \sum_s \chi(-m \cdot k+sx \cdot m) E(x)$$
$$=q^{-d}q^{-1} \sum_{m,x} \sum_{s \not=0} \chi(-m \cdot k+sx \cdot m)
E(x)$$
$$=q^{-1} \sum_{s \not=0} E(s^{-1}k)$$
$$=q^{-1} \sum_{s \not=0} E(sk)=q^{-1}|E \cap l_k|,$$ if $k \not=(0, \dots,
0)$ and
$$ q^{-1}|E|,$$ if $k=(0, \dots, 0)$.

Rewriting (\ref{stage2}) and applying the Parseval identity we get
$$ |E| \sum_m F(m)G(m)=|E|q^d \sum_k \widehat{F}(k)
\overline{\widehat{G}(k)}$$
$$=|E|q^{2d-1} \sum_{k \neq(0, \dots, 0)} |E \cap l_k|
{|\widehat{E}(k)|}^2+|E|q^{2d-1} \cdot |E| \cdot q^{-2d}{|E|}^2$$
$$ \leq C_{geom}|E|q^{2d-1} q^{\frac{\alpha}{d}}
q^{-d}|E|+{|E|}^4q^{-1}C_{geom}{|E|}^2q^{d-1+\frac{\alpha}{d}}+{|E|}^4q^{-1}.$$

Since
$$ {|E|}^4={\left( \sum_t \nu(t) \right)}^2 \leq |\{x \cdot y: x,y \in E\}|
\cdot \sum_t \nu^2(t)$$
$$ \leq |\{x \cdot y: x,y \in E\}|
\left(C_{geom}{|E|}^2q^{d-1+\frac{\alpha}{d}}+
{|E|}^4q^{-1}\right),$$ it follows that
\begin{equation} \label{keylowerbound} |\{x \cdot y: x,y \in E\}| \ge
q \cdot \frac{{|E|}^2}{C_{geom}q^{d+\frac{\alpha}{d}}+{|E|}^2}=r_q \cdot
q.\end{equation}

Suppose that
$$ |E| \ge C_{size}q^{\frac{d}{2}+\frac{\alpha}{2d}}.$$

It follows that
$$ r_q \ge \frac{C_{size}^2}{C_{size}^2+C_{geom}},$$ as desired.

\vskip.125in

\section{Proof of the main arithmetic result (Theorem \ref{kickass})}

\vskip.125in

Let $E=A \times A \times \dots \times A$. The proof of the first part of
Theorem \ref{kickass} follows instantly. To prove the second part observe
that
$$ |E \cap l_y| \leq |A|={|E|}^{\frac{1}{d}}$$ for every $y \in E$.

Then the line (\ref{keylowerbound}) takes the form
$$ |\{(x \cdot y: x,y \in E\}| \ge q \cdot \frac{{|E|}^2}{q^d \cdot
{|E|}^{\frac{1}{d}}+{|E|}^2}.$$

The proof of Theorem \ref{main+} tells us at this point that
$$|\{x \cdot y: x,y \in E\}| \ge q \cdot
\frac{C^{2-\frac{1}{d}}_{size}}{C^{2-\frac{1}{d}}_{size}+1} $$
if
$$|E| \ge C_{size}q^{\frac{d}{2}+\frac{d}{2(2d-1)}}.$$

It follows that if
$$ |A| \ge C^{\frac{1}{d}}_{size}q^{\frac{1}{2}+\frac{1}{2(2d-1)}},$$ then
$$ |dA^2| \ge q \cdot
\frac{C^{2-\frac{1}{d}}_{size}}{C^{2-\frac{1}{d}}_{size}+1}$$ as
desired. This
completes the proof of Theorem \ref{kickass}.

\end{document}